%histabrut.tex:
%%a Plain TeX file by Doron Zeilberger (3 pages)

%begin macros

\baselineskip=14pt
\parskip=10pt

\font\eightrm=cmr8 
\font\eighttt=cmtt8
\magnification=\magstephalf

\def\1{{\overline{1}}}
\def\2{{\overline{2}}}
\parindent=0pt
\overfullrule=0in
\def\Tilde{\char126\relax}
\def\frac#1#2{{#1 \over #2}}
%\headline={\rm  \ifodd\pageno  \RightHead  \else  \LeftHead  \fi}
%\def\RightHead{\centerline{
%Title
%}}
%\def\LeftHead{ \centerline{Doron Zeilberger}}
%end macros
\bf
\centerline
{
HISTABRUT: A Maple Package for  Symbol-Crunching in Probability theory}
\rm
\bigskip
\centerline{ {\it
Doron 
ZEILBERGER}\footnote{$^1$}
{\eightrm  \raggedright
Department of Mathematics, Rutgers University (New Brunswick),
Hill Center-Busch Campus, 110 Frelinghuysen Rd., Piscataway,
NJ 08854-8019, USA.
%\break
{\eighttt zeilberg  at math dot rutgers dot edu} ,
\hfill \break
{\eighttt http://www.math.rutgers.edu/\~{}zeilberg/} .
First version: Aug. 25, 2010.
This article  accompanies the Maple package
{\eighttt http://www.math.rutgers.edu/\Tilde zeilberg/tokhniot/HISTABRUT} .
Exclusively published in the Personal Journal of Shalosh B. Ekhad and Doron Zeilberger
({\eighttt http://www.math.rutgers.edu/\Tilde zeilberg/pj.html}) and {\eighttt http://arxiv.org}.
Supported in part by the NSF.
}
}

{\bf Abstract:} A Maple package {\tt HISTABRUT} (available from
{\eighttt http://www.math.rutgers.edu/\Tilde zeilberg/tokhniot/HISTABRUT})
is presented and briefly described. 
It uses the polynomial ansatz to
discover (often fully rigorously, but in some cases only semi-rigorously (yet rigorizably!)) explicit 
asymptotic formulas for the moments of uni-variate and, more impressively, bi-variate,
discrete probability random variables. It would be hopefully extended, in the future,
to {\it multi-variate} random variables.

Many sequences of discrete random variables (e.g. tossing a (fair or loaded) coin $n$ times, and
keeping track of the number of Heads minus the number of Tails)
are {\it asymptotically normal}. In [Z1], I introduced and described
Maple packages, {\tt CLT} and {\tt AsymptoticMoments}, that empirically-yet-rigorously prove asymptotic normality
for a wide class of sequences of discrete random variables. They used the method of moments.
Furthermore, they are able to prove much stronger theorems than mere ``asymptotic normality'' by
finding the asymptotics (to any desired order!) of the (normalized) moments, rather than only the
leading terms (that should be those of the normal distribution $e^{-x^2/2}/\sqrt{2 \pi}$,
namely $1 \cdot 3 \cdot 5 \cdots (2r-1)$ for the even $2r$-th moment, and
$0$ for the odd moments).

But not {\it all} discrete probability random variables are asymptotically normal!
For example, the number of times that your current capital is positive, upon tossing
a fair coin $n$ times and winning a dollar if it is Heads and losing a dollar it is
Tails, that converges to Paul L\'evy's {\it arcsine distribution} (see [Z2]), and
the other random variables considered by Feller (see [Z3]).
Another intriguing random variable is the {\it duration} of a {\it gambler's ruin}
considered in [Z4].

The much larger Maple package {\tt HISTABRUT} available from

{\eighttt http://www.math.rutgers.edu/\Tilde zeilberg/tokhniot/HISTABRUT}

can handle {\it any} sequence of discrete probability distributions, that the users have to program themselves.
There are quite a few ones pre-programmed, (type {\tt EzraPGF(); } in the Maple package
{\tt HISTABRUT} for a list). It can also sketch the limiting distributions, using procedure
{\tt plotDist} (see the on-line help).

Another new feature is that it can handle {\it directly} sequences of probability distributions
defined in terms of rational generating functions, $R(t,s)$, where the coefficient of $s^n$
in the power-series expansion of $R(t,s)$ in terms of $s$ is the probability generating function
(in $t$) for a typical memebr of a sequence of random variables parametrized by $n$. For example, for tossing a fair coin $n$ times
$R(t,s)=1/(1-s(t+1/t)/2)$. 
Recall that the  Goulden-Jackson[GJ] method (beautifully exposited and extended in [NZ]), and also
included in {\tt HISTABRUT},
outputs such rational functions for the random variable ``number of occurrences of a prescribed  (consecutive) subword''.
First {\tt HISTABRUT} quickly and effortlessly computes {\it explicit} (symbolic) expressions for the 
mean and variance. This is no big deal, and even {\it you}, my dear human readers, can probably do it in 
many cases. Having done this easy task, {\tt HISTABRUT} goes on and computes
the (normalized) even and odd ($2r$-th and $(2r+1)$-th respectively), to {\it any} desired order,
as expressions in {\bf both} $n$ and $r$. Now this is really impressive, and a triumph to experimental mathematics.
It first ``just'' guesses such expressions, but {\it a posteriori}, just by (fully rigorous!) ``hand-waving''
justifies its guesses, by saying that checking a certain number of special cases suffices to prove the
conjectured explicit formulas rigoroulsy. The justification is that at the {\it end of the day}, everything
boils down to {\it polynomial identities}, and we all know that two polynomials of degree $\leq d$ are
identically equal if they coincide in $d+1$ different values. 
In particular it, in any given case, {\it rigorously} reproves the
well-known fact that the distribution is asymptotically normal, but {\it in addition} supplies much more
information, by outputting higher-order asymptotics.

But the most salient new feature is the handling of sequences of {\it bi-variate} discrete random variables, for example
the number of occurrences of two  different words as (consecutive) subwords.
Here it only gives polynomial expressions, in $n$, for the $(r,s)$-mixed moments, for $r,s \leq R$, and
$R$ is a numeric positive integer inputted by the user, but is unable (yet) to find general expressions
in terms of $r$ and $s$.
Here, too, the Goulden-Jackson method, that is built-in, supplies lots of examples.

Whenever the sequence of bivariate discrete probability distributions
happens to be asymptoically {\it independently} normal,  procedure {\tt AnalyseMoms2} can
find explicit expressions, in $n$, of course, but also in {\it both} $r$, and $s$,
for the asymptotic order, to any desired order, for the normalized mixed $(r,s)$ moments.
[More precisely, it finds four distinct expressions for the $(2r,2s)$, $(2r+1,2s)$, $(2r,2s+1)$, and $(2r+1,2s+1)$
mixed moments.]

{\bf Sample input and output}

The ``front'' of the present article

{\tt http://www.math.rutgers.edu/\Tilde zeilberg/mamarim/mamarimhtml/histabrut.html}  

has numerous sample input and output files. The readers  are welcome to edit the input files
in order to produce their own output.

{\bf Future Directions}

Procedure {\tt AnalyseMoms2} (and the verbose version {\tt AnalyseMoms2V})
can only handle sequences of bivariate discrete distributions that are
asymptotically {\it independently} normal. It would be nice to extend it to pairs of
random variables that are non-independently asymptotically normal. This would
first require finding the asymptotic correlation (already done!), and then finding expressions
for the mixed moments for the limiting continuous bi-variate distributions $exp(-x^2/2-y^2/2+bxy)$
where $c=b/(1-b^2)$ is the limiting correlation coefficient.
These are all things that I know how to teach the computer how to do, but I currently don't have time.

Another worthwhile extension is to consider {\it tri}-variate, {\it quad}-variate, and in general, {\it multi}-variate 
sequences of discrete probability distributions.

{\bf References}

[GJ] Ian Goulden  and David M. Jackson,
{\it An inversion theorem for cluster decompositions of
sequences with distinguished subsequences},
J. London Math. Soc.(2){\bf 20} (1979), 567-576.

[NZ] John Noonan and Doron Zeilberger, 
{\it The Goulden-Jackson cluster method: extensions, applications, and implementations},
J. Difference Eq. Appl. {\bf 5} (1999), 355-377. \hfill\break
Available from {\tt http://www.math.rutgers.edu/\Tilde zeilberg/mamarim/mamarimhtml/gj.html} .

[Z1] Doron Zeilberger,
{\it 
The Automatic Central Limit Theorems Generator (and Much More!)},
``Advances in Combinatorial Mathematics: Proceedings of the Waterloo Workshop in Computer Algebra 2008 in honor of Georgy P. Egorychev'', 
chapter 8, 165-174, (I.Kotsireas, E.Zima, eds. Springer Verlag, 2009.). \hfill\break
Available from {\tt http://www.math.rutgers.edu/\Tilde zeilberg/mamarim/mamarimhtml/georgy.html} .

[Z2] Doron Zeilberger,  {\it A new proof that there are $2^n$ ways to toss a coin n times},
Personal Journal of Shalosh B. Ekhad and Doron Zeilberger, \hfill \break
{\tt http://www.math.rutgers.edu/\Tilde zeilberg/mamarim/mamarimhtml/toss.html} .

[Z3] Doron Zeilberger, {\it Fully AUTOMATED Computerized Redux of Feller's (v.1) Ch. III (and Much More!)},
Personal Journal of Shalosh B. Ekhad and Doron Zeilberger, \hfill \break
{\tt http://www.math.rutgers.edu/\Tilde zeilberg/mamarim/mamarimhtml/feller.html} .

[Z4] Doron Zeilberger,
{\it
Symbol Crunching with the Gambler's Ruin Problem},
``Tapas in Experimental Mathematics'' (Tewodros Amdeberhan and Victor Moll, eds.), Contemporary Mathematics, {\bf 457} (2008), 285-292. \hfill\break
Available from {\tt http://www.math.rutgers.edu/\Tilde zeilberg/mamarim/mamarimhtml/ruin.html} .

\end